\documentclass[a4paper,11pt,twoside]{amsart}
\usepackage[english]{babel}
\usepackage[utf8]{inputenc}

\usepackage[a4paper,inner=3cm,outer=3cm,top=4cm,bottom=4cm,pdftex]{geometry}
\usepackage{fancyhdr}
\usepackage{color}
\usepackage{bold-extra}
\usepackage{ mathrsfs }

\usepackage{comment}
\usepackage{graphics}
\usepackage{aliascnt}
\usepackage[pdftex,citecolor=green,linkcolor=black,linkbordercolor={1 1 1}]{hyperref}

\usepackage{amsmath}
\usepackage{amsfonts}
\usepackage{amssymb}
\usepackage{amsthm}
\usepackage{comment}
\usepackage{mathtools}
\usepackage{enumerate}
\usepackage{enumitem} 
\usepackage{bm}

\newtheorem{theorem}{Theorem}[section]

\newtheorem{lemma}[theorem]{Lemma}

\theoremstyle{definition}

\newtheorem{question}[theorem]{Question}
\numberwithin{equation}{section}

\newcommand\eps{\varepsilon}

\newcommand\F{\mathbb{F}}
\newcommand\E{\mathbb{E}}

\newcommand\R{\mathbb{R}}
\newcommand\Z{\mathbb{Z}}

\newcommand\C{\mathbb{C}}

\newcommand\T{\mathbb{T}}

\newcommand{\CC}{\mathcal{C}}

\newcommand{\CL}{\mathcal{L}}

\newcommand{\wh}{\widehat}

\newcommand{\ve}{\mathbf}

\begin{document}
\title{Gowers-uniform sets with few arithmetic progressions}

\author[Shao]{Xuancheng Shao}
\address{Department of Mathematics, University of Kentucky\\
Lexington, KY 40506\\
USA}
\email{xuancheng.shao@uky.edu}
\thanks{XS was supported by NSF grant DMS-2452462.}

\begin{abstract}
For each fixed $s\geq 2$, we construct sparse sets $A\subset \Z/N\Z$ such that the normalized indicator function $\nu_A$ of $A$ satisfies $\|\nu_A-1\|_{U^s}=o(1)$, while the number of three-term arithmetic progressions in $A$ is $o(|A|^3/N)$. This answers the counting form of a question of Gowers and Green, as formulated by Conlon, Fox and Zhao, on the pseudorandomness hypotheses needed in a relative Szemer\'{e}di theorem.
\end{abstract}

\maketitle

\section{Introduction}

Szemer\'edi's theorem \cite{Sze} states that every set of integers of positive upper density contains nontrivial $k$-term arithmetic progressions ($k$-APs), for any fixed $k \geq 3$. By the averaging argument of Varnavides \cite{Var}, a set of density at least $\alpha>0$ in $\Z/N\Z$ contains at least $c(k,\alpha)N^2$ $k$-APs for all sufficiently large primes $N$. Gowers's proof \cite{Gow} gives a criterion for the count to be close to that of a random set. Indeed, if $A\subset\Z/N\Z$ has density $\alpha$, the generalized von Neumann theorem (see, for example, \cite[Section~7]{CFZ}) implies that
\[
\big|\Lambda_k(1_A)-\alpha^k\big|
   \leq k\|1_A-\alpha\|_{U^{k-1}}.
\]
Thus $U^{k-1}$-uniformity (i.e. $\|1_A-\alpha\|_{U^{k-1}} = o(1)$) determines the asymptotic count for dense sets. Here the $k$-AP counting operator $\Lambda_k$ is defined by
$$
\Lambda_k(f) = \E_{x,d \in G}\prod_{j=0}^{k-1}f(x+jd)
$$
for a function $f:G\rightarrow\C$ on a finite abelian group $G$. For $s \geq 2$, the Gowers $U^s$-norm of $f$ is defined by
$$
\|f\|_{U^s}^{2^s} = \E_{x\in G}\E_{\ve{h}\in G^s} \prod_{\bm{\omega}\in\{0,1\}^s} \CC^{|\bm{\omega}|} f(x+\bm{\omega}\cdot\ve{h}),
$$
where $\CC$ denotes complex conjugation and $|\bm{\omega}| = \omega_1+\cdots+\omega_s$ for $\bm{\omega}=(\omega_1,\cdots,\omega_s)$. See \cite{GT} for basic properties of Gowers norms and the central roles they play in additive combinatorics.

In the sparse setting when the density $\alpha$ of $A$ tends to $0$, we will find it convenient to use the normalized indicator function $\nu_A$ of $A$, defined to be $\nu_A=\alpha^{-1}1_A$. Thus $\nu_A$ has average $1$, and the generalized von Neumann theorem gives
\begin{equation}\label{eq:gvn}
\big|\Lambda_k(\nu_A) - 1\big| \leq k\alpha^{-(k-1)} \|\nu_A-1\|_{U^{k-1}}.
\end{equation}
The relative Szemer\'edi theorem of Green and Tao \cite{GT} gives a lower bound for the $k$-AP count in any relatively dense subset of $A$, assuming that $\nu_A$ satisfies some pseudorandomness hypotheses. Conlon, Fox and Zhao \cite{CFZ} showed that a family of linear forms conditions suffices for the pseudorandomness hypothesis. These linear forms conditions essentially include the Gowers uniformity assumptions $\|\nu_A-1\|_{U^{k-1}} = o(1)$, but they also include additional conditions. On the other hand, their companion note \cite{CFZnorm} proves a relative Szemer\'edi theorem where the pseudorandomness hypothesis is $\|\nu_A-1\|_{U^{k-1}} = o(\alpha^{k-1})$ (which immediately implies that $A$ itself has the expected number of $k$-APs by \eqref{eq:gvn}).

For $3$-APs, Green's transference argument for the primes \cite{GreenRoth} and the restriction theorem of Green and Tao \cite{GTRest} give a relative Roth theorem under Fourier uniformity of the majorant and a restriction estimate for the majorant. (See Prendiville's note \cite{Prendiville} on different versions of transference principles.)

Following a question of Gowers and Green (see \cite{GowTr} and \cite[Section~7]{CFZ}), one may ask whether the relative Szemer\'edi theorem holds in $A$, under only the  Gowers uniformity assumption that $\|\nu_A-1\|_{U^s} = o(1)$ for some $s = s(k)$. We show that the answer is negative already for $k=3$, by constructing sparse subsets $A \subset \Z/N\Z$ which are $U^s$-uniform for any fixed $s$, while having few $3$-APs.

\begin{theorem}\label{thm:main-cyclic}
Fix a positive integer $s \geq 2$.
Let $G = \Z/N\Z$ for some prime $N$ sufficiently large in terms of $s$. For any $\eps \in ((\log N)^{-1/2},1/2)$, there exists a subset $A \subset G$ such that its normalized indicator function $\nu_A=(|G|/|A|)1_A$  satisfies
$$
\|\nu_A-1\|_{U^s}^{2^s} = O_s(\eps), \ \ \Lambda_3(\nu_A) \leq 4\exp(-c (\eps \log N)^{1/2})
$$
for some constant $c = c(s) > 0$.
\end{theorem}

In a related direction, Dodos and Kanellopoulos \cite{DK} show that if $\|\nu-1\|_{U^{2s}}=o(1)$, then every $0\leq f\leq\nu$ admits a dense model $0\leq g\leq1$ with $\|f-g\|_{U^s}=o(1)$. On the other hand, the question considered here concerns the progression count.

One can also ask whether $U^s$-uniformity forces the existence of even one nontrivial $3$-AP, which is not addressed by our construction. Bennett and Bohman \cite{BB} give a negative answer for $s=2$. We return to this existence question in Section~\ref{sec:related}.

For longer progressions there are related examples at fixed density. Gowers \cite{GowFour} constructed Fourier-uniform sets of density $1/2+o(1)$ with $4$-AP density at most $1/16-c$ for an absolute constant $c>0$. Wolf \cite{Wolf} obtained a larger deficit. For odd $k\geq5$, Deng, Tidor and Zhao \cite{DTZ} constructed $U^{k-2}$-uniform sets of density $\alpha+o(1)$ with $k$-AP density at most $\alpha^{c_k\log(1/\alpha)}+o(1)$ for sufficiently small fixed $\alpha>0$.

\subsection*{Proof ideas}

Our construction in Theorem \ref{thm:main-cyclic} uses the following simple observation. If $A$ is a subset of a finite abelian group $H$ (equipped with the probability measure) or the circle $\T=\R/\Z$ (equipped with the Haar measure), formed by deleting a suitable set of small measure $\eta$, then three distinct vertices of a random cube
$$
\{x + \bm{\omega}\cdot\ve{h}: \bm{\omega}\in\{0,1\}^s\}
$$
hit the deleted set independently. On the other hand, the three terms of a $3$-AP can all hit the deleted set with probability of order $\eta^2$.
Consequently the deletion changes the (normalized) cube averages by only $O_s(\eta^3)$, while it changes the (normalized) $3$-AP counts by $\gg \eta^2$. This observation is related to the Sidorenko property of linear equations; see \cite{SW,FPZ}. The details of this observation are given in Section \ref{sec:TD}.

Taking an $m$-fold product for some $m$ satisfying $m\eta^3=o(1)$ and $m\eta^2\to\infty$, the cube averages are still $1+o(1)$, while the $3$-AP counts tend to $0$. In $\Z/N\Z$, this tensoring argument is implemented by modeling $\Z/N\Z$ by $\T^m$ via the map $n\mapsto (b_1n/N, \cdots, b_mn/N)$ for appropriate choices of $b_1,\cdots,b_m\in \Z/N\Z$. A similar modeling argument was used by Green and Wolf \cite{GW} in the construction of a $3$-AP-free subset. 
The details of this tensoring argument are given in Section \ref{sec:proof}.

\section{Estimates for a punctured circle}\label{sec:TD}

In this section, we estimate the cube averages and the three-term progression count in a punctured circle $\T\setminus D$, where $\T = \R/\Z$ and $D \subset \T$ is a subset of small measure.

\begin{lemma}\label{lem:torus}
Let $D \subset \T$ be a subset of measure  $\eta \in (0,1/4)$, and let $w: \T\rightarrow\R$ be the normalized indicator function of $\T\setminus D$ defined by $w= (1-\eta)^{-1}1_{\T\setminus D}$.
Then for every fixed $s \geq 2$ and $B \subset \{0,1\}^s$, we have
$$
\int_{x \in \T} \int_{\ve{h}\in\T^s} \prod_{\bm{\omega}\in B}w(x+\bm{\omega}\cdot\ve{h}) = 1 + O_s(\eta^3).
$$
\end{lemma}

\begin{proof}
Since $w = (1-\eta)^{-1}(1-1_D)$, we can write the integral on the left-hand side as
$$
 (1-\eta)^{-|B|} \sum_{B'\subset B} (-1)^{|B'|}\int_{x \in \T} \int_{\ve{h}\in\T^s} \prod_{\bm{\omega}\in B'}1_D(x+\bm{\omega}\cdot\ve{h}).
$$
We claim that
\begin{equation}\label{eq:torus-claim}
\int_{x \in \T} \int_{\ve{h}\in\T^s} \prod_{\bm{\omega}\in B'}1_D(x+\bm{\omega}\cdot\ve{h})=\begin{cases} \eta^{|B'|}, & \text{if }|B'| \leq 2, \\ O(\eta^3), & \text{if }|B'| \geq 3.\end{cases}
\end{equation}
For $|B'| \geq 3$, pick three distinct $\bm{\omega}_1,\bm{\omega}_2,\bm{\omega}_3\in B'$. Since $\bm{\omega}_1,\bm{\omega}_2,\bm{\omega}_3 \in \{0,1\}^s$, the vectors $(1,\bm{\omega}_i)$ for $1\leq i\leq 3$ are linearly independent. Hence the map
$$
(x,\ve{h})\mapsto (x+\bm{\omega}_1\cdot\ve{h},x+\bm{\omega}_2\cdot\ve{h},x+\bm{\omega}_3\cdot\ve{h})
$$
is a surjective linear map from $\T^{s+1}$ to $\T^3$ and takes the Haar measure on $\T^{s+1}$ to the Haar measure on $\T^3$. Thus the simultaneous conditions $x+\bm{\omega}_i\cdot\ve{h}\in D$ for each $1 \leq i \leq 3$ happen with probability $\eta^3$. This proves the claim \eqref{eq:torus-claim} when $|B'| \geq 3$. The claim for $|B'| \leq 2$ follows from the same argument above, using the linear independence of $(1,\bm{\omega})$ for $\bm{\omega}\in B'$.
It follows that 
$$
\int_{x \in \T} \int_{\ve{h}\in\T^s} \prod_{\bm{\omega}\in B}w(x+\bm{\omega}\cdot\ve{h}) = (1-\eta)^{-|B|}\Big(
1-|B|\eta+\binom{|B|}{2}\eta^2+O_s(\eta^3)\Big).
$$
The desired conclusion follows by noting that
$$
1-|B|\eta+\binom{|B|}{2}\eta^2 = (1-\eta)^{|B|} + O_s(\eta^3)
$$
from Taylor expansion. This completes the proof.
\end{proof}

\begin{lemma}\label{lem:torus-3ap}
Let $\eta \in (0,1/4)$, and let $w: \T\rightarrow\R$ be the normalized indicator function of $\T\setminus [0,\eta)$ defined by $w = (1-\eta)^{-1}1_{[\eta,1)}$.
Then we have
$$
\int_{x,d\in\T} w(x) w(x+d) w(x+2d) \leq 1-\frac{\eta^2}{4}.
$$
\end{lemma}

\begin{proof}
Since $w = (1-\eta)^{-1}(1-1_{[0,\eta)})$, we can write the integral on the left-hand side as
$$
(1-\eta)^{-3}\int_{x,d\in\T}(1-1_{[0,\eta)}(x)) (1-1_{[0,\eta)}(x+d)) (1-1_{[0,\eta)}(x+2d)).
$$
Each of $x,x+d,x+2d$ lies in $[0,\eta)$ with probability $\eta$, and each pair lies in $[0,\eta)^2$ with probability $\eta^2$. If all three terms $x, x+d,x+2d$ lie in $[0,\eta)$, then we write $y = x+d$ and so $x+2d=2y-x$. For $x,y \in [0,\eta)$, we have $-\eta \mathrel{<} 2y-x \mathrel{<} 2\eta$.  Since $\eta<1/4$, the condition $2y-x\in[0,\eta)$ modulo $1$ is equivalent to
$$
\frac{x}{2}\leq y\mathrel{<}\frac{x+\eta}{2}.
$$
Thus all three terms lie in $[0,\eta)$ with probability $\eta^2/2$. It follows that
$$
\begin{aligned}
\int_{x,d\in\T} w(x)w(x+d)w(x+2d)
&= (1-\eta)^{-3}(1-3\eta+3\eta^2-\eta^2/2)\\
&= 1-\frac{\eta^2(1/2-\eta)}{(1-\eta)^3}
\leq 1-\frac{\eta^2}{4},
\end{aligned}
$$
since $\eta<1/4$. This completes the proof.
\end{proof}

\section{Proof of Theorem \ref{thm:main-cyclic}}\label{sec:proof}

Let $G = \Z/N\Z$ for some prime $N$ sufficiently large in terms of $s$.
Choose parameters $\eta,m$ by
$$
\eta = \Big(\frac{\eps}{\log N}\Big)^{1/2}, \ \ m = \lfloor 10^{-s}\eta^{-1}\log N\rfloor.
$$
Let $I$ be the interval $\{\lceil\eta N\rceil,\ldots,N-1\}$ in $G$, and let $\nu_I = (1-\eta)^{-1}1_I$. Then $\nu_I(n) = w(n/N)$, where $w$ is the normalized indicator function of $[\eta,1) \subset \T$. For $\ve{b}=(b_1,\cdots,b_m)\in G^m$, define the subset $A_{\ve{b}} \subset G$ by
$$
A_{\ve{b}} = \{n \in G: b_in\in I\text{ for each }1 \leq i \leq m\}.
$$
We will choose $b_1,\cdots,b_m$ independently and uniformly at random. Let $f_{\ve{b}} = (1-\eta)^{-m}1_{A_{\ve{b}}}$. This is not quite the normalized indicator function of $A_{\ve{b}}$, but we will return to it at the end of the argument. Note that
\begin{equation}\label{eq:f-w}
f_{\ve{b}}(n) = \prod_{i=1}^m\nu_I(b_in) = \prod_{i=1}^m w\Big(\frac{b_in}{N}\Big).
\end{equation}

The following lemma allows us to approximate a counting estimate in $G$ by a corresponding one in $\T$.

\begin{lemma}\label{lem:modeling}
Let the notation and assumptions be as above.
Let $\CL$ be a collection of $t \leq 2^s$ nonzero homogeneous linear forms with integral coefficients in $r$ variables. Let
$$
F_{\CL}(\ve{n})=\prod_{L\in\CL} \nu_I(L(\ve{n})),\ \  p_{\CL}=\E_{\ve{n}\in G^r}F_{\CL}(\ve{n}).
$$
Then we have
$$
\E_{\ve{n}\in G^r}\Big|(\E_{b\in G}F_{\CL}(b\ve{n}))^m-p_{\CL}^m\Big|
\ll_s N^{-1/2}
$$
Moreover, if the sum of the absolute values of all the coefficients in each linear form in $\CL$ is $O_s(1)$, then
$$
p_{\CL}^m = \Big(\int_{\ve{x}\in\T^r} \prod_{L\in\CL} w(L(\ve{x})) \Big)^m + O_s(N^{-1/2}).
$$
\end{lemma}

\begin{proof}
First we prove the conclusion for $m=1$ (with better error terms).
We use the Fourier expansion
$$
\nu_I(x) = \sum_{a \in G} \wh{\nu_I}(a) e_N(ax),
$$
where $e_N(t) = \exp(2\pi it/N)$.
Let $\CL=(L_1,\cdots,L_t)$. Then for each $\ve{n}\in G^r$ we have
$$
\begin{aligned}
\E_{b\in G}F_{\CL}(b\ve{n})
&=\sum_{a_1,\cdots,a_t\in G} \Big(\prod_{j=1}^t\wh{\nu_I}(a_j)\Big) \E_{b \in G} e_N\Big(b(a_1L_1(\ve{n})+\cdots+a_tL_t(\ve{n}))\Big) \\
&=\sum_{\substack{a_1,\cdots,a_t\in G\\ a_1L_1(\ve{n})+\cdots+a_tL_t(\ve{n})=0}}\Big(\prod_{j=1}^t\wh{\nu_I}(a_j)\Big).
\end{aligned}
$$
On the other hand, we have
$$
\begin{aligned}
p_{\CL} &= \sum_{a_1,\cdots,a_t\in G} \Big(\prod_{j=1}^t \wh{\nu_I}(a_j)\Big) \E_{\ve{n} \in G^r} e_N(a_1L_1(\ve{n}) + \cdots+a_tL_t(\ve{n})) \\
&= \sum_{\substack{a_1,\cdots,a_t\in G\\ a_1L_1+\cdots+a_tL_t=0}} \Big(\prod_{j=1}^t \wh{\nu_I}(a_j)\Big).
\end{aligned}
$$
Subtracting the two expressions above and averaging over $\ve{n}$, we have
$$
\frac{1}{N^r}\sum_{\ve{n}\in G^r}\Big|\E_{b\in G}F_{\CL}(b\ve{n})-p_{\CL}\Big| \leq \frac{1}{N^r}\sum_{\substack{\ve{n}\in G^r, a_1,\cdots,a_t\in G\\ a_1L_1+\cdots+a_tL_t\neq 0\\ a_1L_1(\ve{n}) + \cdots+a_tL_t(\ve{n})=0}} \prod_{j=1}^t |\wh{\nu_I}(a_j)|.
$$
Given $a_1,\cdots,a_t\in G$ with $a_1L_1+\cdots+a_tL_t\neq 0$, the condition $a_1L_1(\ve{n}) +\cdots+a_tL_t(\ve{n})=0$ is satisfied for exactly a $1/N$ proportion of $\ve{n}\in G^r$. It follows that the right-hand side above is at most
$$
\frac{1}{N} \sum_{a_1,\cdots,a_t\in G} \prod_{j=1}^t |\wh{\nu_I}(a_j)| = \frac{(1-\eta)^{-t}}{N}\Big(\sum_{a\in G}|\wh{1_I}(a)|\Big)^t \ll_t \frac{(\log N)^t}{N}
$$
since $I$ is an interval. This proves the first conclusion when $m=1$:
$$
\E_{\ve{n}\in G^r}\Big|\E_{b\in G}F_{\CL}(b\ve{n})-p_{\CL}\Big|
\ll_s \frac{(\log N)^{O_s(1)}}{N}.
$$

For the second conclusion, divide $\T^r$ into $N^r$ boxes of side lengths $1/N$ of the form $\ve{n}/N + [0,1/N)^r$, where $\ve{n} \in G^r$. On each such box, we have
$$
\prod_{L\in\CL} w(L(\ve{x})) = \prod_{L\in\CL} \nu_I(L(\ve{n}))
$$
for all $\ve{x}$ in the box, unless $L(\ve{n}/N)$ is within $O_s(1/N)$ of the boundary of $[\eta,1)$ for some $L \in \CL$. Since this happens for $O_s(1/N)$ proportion of $\ve{n}$, we have
$$
\int_{\ve{x}\in\T^r} \prod_{L\in\CL} w(L(\ve{x})) = \E_{\ve{n}\in G^r} \prod_{L\in\CL} \nu_I(L(\ve{n})) + O_s\Big(\frac{1}{N}\Big) = p_{\CL} + O_s\Big(\frac{1}{N}\Big).
$$
This proves the second conclusion when $m=1$.

For the value of $m$ chosen at the beginning of this section, we use the inequality
$$
|u^m-v^m| \leq m|u-v|\cdot \max(|u|^{m-1},|v|^{m-1})
$$
with the two choices 
$$
(u,v) = \Big(p_{\CL}, \E_{b \in G}F_{\CL}(b\ve{n})\Big) \text{ and }\Big(p_{\CL}, \int_{\ve{x}\in\T^r}\prod_{L\in\CL}w(L(\ve{x}))\Big).
$$
For both choices, we have $\max(|u|,|v|) \leq (1-\eta)^{-t}$. Since
$$
m(1-\eta)^{-tm} \leq m(1-\eta)^{-2^sm} \leq N^{0.1}
$$
by our choice of $m$, both conclusions hold.
\end{proof}

We now apply Lemma~\ref{lem:modeling} to estimate the density of $A_{\ve{b}}$, the $U^s$-norm of $f_{\ve{b}}-1$, and the $3$-AP count. Lemma \ref{lem:modeling} reduces these counting estimates in $G$ to a corresponding counting integral in $\T^m$.

\subsection*{The density of $A_{\ve{b}}$}

First we prove that
\begin{equation}\label{eq:proof-density}
\E f_{\ve{b}}  =(1-\eta)^{-m} \frac{|A_{\ve{b}}|}{N} = 1 + O_s(N^{-1/4}).
\end{equation}
with probability at least $3/4$. By first using \eqref{eq:f-w} and then applying Lemma \ref{lem:modeling} with the single linear form $L(n)=n$, we have
$$
\E_{\ve{b}\in G^m}\E_{n \in G} f_{\ve{b}}(n) = \E_{n\in G}\Big(\E_{b\in G} \nu_I(bn)\Big)^m = \Big(\int_{x\in\T} w(x)\Big)^m + O_s(N^{-1/2}) = 1+O_s(N^{-1/2}).
$$
For the second moment, by \eqref{eq:f-w} we have
$$
\E_{\ve{b}\in G^m} \Big(\E_{n\in G}f_{\ve{b}}(n)\Big)^2 = \E_{n_1,n_2\in G}\Big(\E_{b \in G} \nu_I(bn_1)\nu_I(bn_2)\Big)^m.
$$
By Lemma \ref{lem:modeling} applied to the two linear forms $(n_1,n_2)\mapsto n_1$ and $(n_1,n_2)\mapsto n_2$, the expression above is
$$
\Big(\int_{x_1,x_2\in\T} w(x_1)w(x_2)\Big)^m + O_s(N^{-1/2}) = 1+O_s(N^{-1/2}).
$$
It follows that
$$
\E_{\ve{b}\in G^m}\Big(\E_{n \in G}f_{\ve{b}}(n)-1\Big)^2 = O_s(N^{-1/2}).
$$
By Markov's inequality, we have
$$
\Big|\E_{n \in G}f_{\ve{b}}(n)-1\Big| = O_s(N^{-1/4})
$$
with probability at least $3/4$. This gives \eqref{eq:proof-density}.

\subsection*{Gowers uniformity of $f_{\ve{b}}-1$}

Next we prove that 
\begin{equation}\label{eq:proof-gowers}
\|f_{\ve{b}}-1\|_{U^s}^{2^s} = O_s(\eps)
\end{equation}
with probability at least $3/4$. By expanding out the product in the definition of the $U^s$-norm of $f_{\ve{b}}-1$, we have
$$
\|f_{\ve{b}}-1\|_{U^s}^{2^s} = \sum_{B \subset \{0,1\}^s} (-1)^{|B|}\cdot  \E_{x\in G}\E_{\ve{h}\in G^{s}} \prod_{\bm{\omega}\in B}f_{\ve{b}}(x+\bm{\omega}\cdot\ve{h}).
$$
Averaging over $\ve{b}$ and using \eqref{eq:f-w}, we have
$$
\E_{\ve{b} \in G^m}\|f_{\ve{b}}-1\|_{U^s}^{2^s} =
\sum_{B \subset \{0,1\}^s} (-1)^{|B|}\cdot  \E_{x\in G}\E_{\ve{h}\in G^{s}} \Big(\E_{b \in G}\prod_{\bm{\omega}\in B}\nu_I(b(x+\bm{\omega}\cdot\ve{h}))\Big)^m.
$$
For each $B\subset\{0,1\}^s$, we apply Lemma \ref{lem:modeling} to the collection of linear forms consisting of $(x,\ve{h})\mapsto x+\bm{\omega}\cdot\ve{h}$ for $\bm{\omega}\in B$. Thus the average over $x$ and $\ve{h}$ above is
$$
\Big(\int_{x\in\T}\int_{\ve{h}\in\T^s} \prod_{\bm{\omega}\in B} w(x+\bm{\omega}\cdot\ve{h}) \Big)^m + O_s(N^{-1/2}).
$$
By Lemma \ref{lem:torus}, the expression above is
$$
(1 + O_s(\eta^3))^m + O_s(N^{-1/2}) = 1 + O_s(m\eta^3 + N^{-1/2}) = 1 + O_s(\eps)
$$
by our choice of $m$ and $\eta$.
It follows that
$$
\E_{\ve{b}\in G^m}\|f_{\ve{b}}-1\|_{U^s}^{2^s} = \sum_{B\subset\{0,1\}^s} (-1)^{|B|}(1+O_s(\eps)) = O_s(\eps).
$$
Thus \eqref{eq:proof-gowers} follows by Markov's inequality.

\subsection*{\texorpdfstring{The $3$-AP count}{The 3-AP count}}

Finally we prove that
\begin{equation}\label{eq:proof-3ap}
\Lambda_3(f_{\ve{b}}) \leq 4\exp(-c(\eps \log N)^{1/2})
\end{equation}
with probability at least $3/4$ for some constant $c=c(s)>0$. By \eqref{eq:f-w}, we have
$$
\E_{\ve{b}\in G^m} \Lambda_3(f_{\ve{b}})= \E_{x,d\in G}\Big(\E_{b\in G} \nu_I(bx)\nu_I(b(x+d))\nu_I(b(x+2d))\Big)^m.
$$
We apply Lemma \ref{lem:modeling} to the collection of the three linear forms $x, x+d, x+2d$ in the variables $x,d$. Thus
$$
\E_{\ve{b}\in G^m} \Lambda_3(f_{\ve{b}}) = \Big(\int_{x,d\in\T} w(x) w(x+d) w(x+2d)\Big)^m + O_s(N^{-1/2}).
$$
By Lemma \ref{lem:torus-3ap}, we have
$$
\E_{\ve{b}\in G^m} \Lambda_3(f_{\ve{b}})  \leq \Big(1-\frac{\eta^2}{4}\Big)^m + O_s(N^{-1/2}) \leq \exp(-c(\eps \log N)^{1/2})
$$
for some constant $c=c(s)>0$, by our choice of $m$ and $\eta$.
Thus \eqref{eq:proof-3ap} follows by Markov's inequality.

\subsection*{Conclusion of proof}

It follows that there exists a choice $\ve{b} \in G^m$ such that \eqref{eq:proof-density}, \eqref{eq:proof-gowers}, and \eqref{eq:proof-3ap} hold simultaneously. Fix one such $\ve{b}$, take $A = A_{\ve{b}}$, and let $\nu_A$ be the normalized indicator function of $A$. Then $\nu_A=f_{\ve{b}}/M$, where $M = \E f_{\ve{b}} = 1+O_s(N^{-1/4})$ by \eqref{eq:proof-density}. Thus $\Lambda_3(\nu_A) = M^{-3}\Lambda_3(f_{\ve{b}})$ and 
$$
\|\nu_A-1\|_{U^s} \leq M^{-1}(\|f_{\ve{b}}-1\|_{U^s}+|M-1|).
$$
Hence the estimates \eqref{eq:proof-gowers} and \eqref{eq:proof-3ap} continue to hold with $f_{\ve{b}}$ replaced by $\nu_A$ (after decreasing the constant $c(s)$ in \eqref{eq:proof-3ap} slightly).

\section{Related questions}\label{sec:related}

In this section, we discuss the existence question in cyclic groups $G = \Z/N\Z$ (with $N$ a large prime) and in the finite field model $G = \F_3^n$ (with $n$ large). Let $N = |G|$ in both settings.

\begin{question}[Existence]
Is there an integer $s\geq 2$ such that the following statement holds? If $A\subset G$ satisfies $\|\nu_A-1\|_{U^s}=o(1)$, then $A$ must contain a nontrivial $3$-AP.
\end{question}

\begin{question}[Relative existence]
Is there an integer $s\geq 2$ such that, for every fixed $\delta>0$, the following statement holds? If $A\subset G$ satisfies $\|\nu_A-1\|_{U^s}=o(1)$, then every subset $A'\subset A$ of relative density at least $\delta$ must contain a nontrivial $3$-AP.
\end{question}

In both groups the answer is known to be negative for $s=2$. In $\Z/N\Z$, Bennett and Bohman \cite{BB} use the random greedy process to construct a Fourier-uniform set $A$ containing no nontrivial $3$-APs, and their method leads similarly to a construction in the finite field model. More generally, for $k=2^{s-1}+1$ with $s\geq2$, their construction gives a $k$-AP-free set satisfying $\|\nu_A-1\|_{U^s}=o(1)$.

In $\F_3^n$, Lev's cap-set construction \cite{Lev}, also described in \cite{EL}, gives a deterministic example. Assume that $n$ is divisible by $3$ and take $q = 3^{n/3}$. In the finite field $\F_q$, choose a nonsquare $D\in\F_q$, and let $A$ be the graph of $Q(u,v)=u^2-Dv^2$:
$$
A = \{(u,v, u^2-Dv^2): u,v \in \F_q\} \subset \F_q^3.
$$
We can view $A$ as a subset of $\F_3^n$ by identifying $\F_q$ with $\F_3^{n/3}$ in the natural way (as vector spaces over $\F_3$). If $A$ contains a $3$-AP of the form
$$
(u-h,v-k,Q(u-h,v-k)), \ \ (u,v,Q(u,v)), \ \ (u+h,v+k,Q(u+h,v+k))
$$
for some $(h,k)\neq (0,0)$, then
$$
2Q(h,k) = Q(u-h,v-k)+Q(u+h,v+k)-2Q(u,v) = 0,
$$
and so $h=k=0$ (since $D$ is nonsquare). Thus $A$ contains no nontrivial $3$-APs.  For the Fourier uniformity, note that $\nu_A = q1_A$ since $|A|=q^2$. By Fourier analysis on $\F_q^3$, we have
$$
\wh{\nu_A}(a,b,c) = \frac{1}{q^2}\sum_{u,v\in\F_q} \psi(au + bv + c(u^2-Dv^2))
$$
for $a,b,c\in\F_q$, where $\psi(t) = e_3(-\operatorname{Tr}_{\F_q/\F_3}(t))$. If $c=0$, then the sum vanishes unless $(a,b)=(0,0)$. If $c\neq 0$, then
$$
\wh{\nu_A}(a,b,c) = \frac{1}{q^2}\Big(\sum_{u \in \F_q} \psi(au +cu^2)\Big) \Big(\sum_{v \in \F_q} \psi(bv -cDv^2)\Big),
$$
and each of the two sums is a quadratic Gauss sum of size $q^{1/2}$. Thus $|\wh{\nu_A}(a,b,c)| = 1/q$ when $c \neq 0$. Consequently
\[
\|\nu_A-1\|_{U^2}^4
=\sum_{(a,b,c)\neq (0,0,0)}|\widehat{\nu_A}(a,b,c)|^4
=q^2(q-1)q^{-4},
\]
which is $o(1)$ as $n\rightarrow\infty$.

Finally, we remark that $s=5$ suffices in the finite field model, even for the relative existence question. Indeed, if $\|\nu_A-1\|_{U^5}=o(1)$ then the density $\alpha$ of $A$ must satisfy
$$
\alpha^{-1} = \|\nu_A^2\|_{U^1} \leq \|\nu_A^2\|_{U^4} \leq N^{1/16}\|\nu_A\|_{U^5}^2 \leq (1+o(1))N^{1/16},
$$
where the middle inequality follows from the recursive definition of the Gowers norms:
$$
\|\nu_A\|_{U^5}^{32} = \E_{h\in G} \|\Delta_h\nu_A\|_{U^4}^{16} \geq \frac{1}{N}\|\Delta_0\nu_A\|_{U^4}^{16} = \frac{1}{N}\|\nu_A^2\|_{U^4}^{16}.
$$
Here $\|f\|_{U^1} = |\E f|$ and $\Delta_hf(x) = f(x+h)\overline{f(x)}$.
Hence $|A| \gg N^{15/16} = 3^{(15/16)n}$, and similarly for any subset of $A$ of fixed positive relative density. This ensures the existence of $3$-APs by the cap-set bound of Ellenberg and Gijswijt \cite{EG}, giving that any $3$-AP-free subset $B\subset \F_3^n$ must satisfy $|B|\ll 2.756^n$ (and $2.756<3^{15/16}$). Ellenberg and Gijswijt's proof builds on the polynomial method of Croot, Lev and Pach \cite{CLP} for progression-free subsets of $(\Z/4\Z)^n$.

\end{document}